\documentclass[12pt]{amsart}
   
\usepackage{a4}

\usepackage{enumerate}

\newcommand{\DB}{\mathbb{B}}
\newcommand{\E}{\mathbb{E}}

\newcommand{\DR}{\mathbb{R}}
\newcommand{\DS}{\mathbb{S}}

\newcommand{\DEins}{\mathbf{1}}

\newcommand{\norm}[1][\cdot]{\vert \hspace{-0.15em} \vert \hspace{-0.1em} #1 \hspace{-0.1em} \vert \hspace{-0.15em} \vert}

\newcommand{\scalp}[1][\cdot \: ,\cdot]{\langle #1 \rangle}
\newcommand{\Scalp}[1][\cdot \: ,\cdot]{\left \langle #1 \right \rangle}

\newcommand{\sqfr}[2]{\sqrt{\frac{#1}{#2}}}

\newcommand{\CA}{\mathcal{A}}

\newcommand{\CL}{\mathcal{L}}

\newcommand{\AS}{\CA}
\newcommand{\LS}{\CL}

\newcommand{\spann}{\mathrm{span}}
\newcommand{\op}[2][C]{#1 | #2}

\DeclareMathOperator{\conv}{conv}

\newtheorem{Satz}{Theorem}[section]
\newtheorem{Proposition}[Satz]{Proposition}
\newtheorem{Lemma}[Satz]{Lemma}
\newtheorem{Folgerung}[Satz]{Corollary}
\newtheorem{Beispiel}[Satz]{Example}

\newtheorem{Definition}[Satz]{Definition}

\title{radii of regular polytopes}
\author[R.~Brandenberg]{Ren\'{e} Brandenberg}
\address{Ren\'{e} Brandenberg, Zentrum Mathematik, Technische
  Universit\"at M\"unchen, D-80290 Munich, Germany} 
\email{brandenb@mathematik.tu-muenchen.de}
\urladdr{http://www-m9.mathematik.tu-muenchen.de/\~{}brandenb/}
\address{currently on leave at Technische Universit\"at Wien \\
        Institut f\"ur Analysis und Technische Mathematik \\
        Wiedner Hauptstr. 8--10 \\
        A--1040 Wien}

\thanks{Research of the author was supported by the European
  Union, through a Marie Curie Fellowship, Contract-No.: HPMT-CT-2000-00037}

\thanks{The author wants to thank Keith Ball for usefull discussions.}


\date{\today}

\begin{document}

\begin{abstract}
This paper deals with the three types of regular polytopes which
exist in all dimensions -- regular simplices, cubes and
regular cross-polytopes -- and their outer and inner radii. While the
inner radii of regular simplices are well studied, only a few cases 
are solved for the outer radii. 
We give a lower bound on these radii, and show that
this bound is tight in almost 3 out of 4 dimensions. In a further
section we complete the results about inner and outer radii of general
boxes and cross-polytopes. Finally, because cubes and regular
cross-polytopes are radii-minimal 
projections of simplices, we show that it is possible to deduce the results
about their radii from the results about the outer radii of simplices.
\end{abstract}

\maketitle

\section{Introduction and basic notations} \label{basics}

 There are three classes of regular polytopes which last in general
 $d$-space: regular simplices, (hyper-) cubes, and regular cross-polytopes. In
 this paper we investigate the inner and outer $j$-radii of this
 polytopes. Here the inner $j$-radius of a body is defined as the
 radius of a biggest $j$-ball fitting into the body, and the outer
 $j$-radius as the smallest radius of a circumball of an orthogonal
 projection of the body onto any $j$-space. 

 While the inner radii of
 regular simplices are well studied \cite{ball-92}, very less is known
 about their outer radii. We give a lower bound on these radii,
 and show that this bound is tight in almost 3 out of 4 cases. An
 important step towards this result is the investigation of quasi 
 isotropic polytopes (Kawashima called them $\pi$-polytopes
 \cite{kawashima-91}, but we prefer to call them isotropic as they are
 in an isotropic position in the sense of
 \cite{giannopoulos-milman-00}). Specifically, we will show that the
 existence of a quasi-isotropic $j$-dimensional polytope with $d+1$
 vertices is 
 equivalent to the existence of a projection of the regular simplex
 such that the lower bound is attained.

 In a further section we investigate the radii of general boxes and
 cross-polytopes. While the inner radii of boxes were
 computed in \cite{esvw-98}, nothing could be found about their
 outer radii in the literature. We close this gap, as
 well as we transfer the results about boxes to results about
 cross-polytopes via polarization. Finally, we show that the radii of
 cubes and regular cross-polytopes can be obtained almost completely 
 from our results about the outer-radii of regular simplices.  

   Let $\E^d= (\DR^d,\norm)$ denote the $d$-dimensional Euclidean
   space, $d \ge 2$, 
   $\DB$ and $\DS$ the unit ball and the unit sphere in $\E^d$, and
   $\scalp$ the usual \emph{scalar product}
   $\scalp[x,y]=x^Ty$. Furthermore, we use 
  $\{e_1,\dots,e_d\}$ for the \emph{standard basis} of $\E^d$. 

   A set
   $C \subset \E^d$ is called a \emph{body} if it is 
   bounded, closed, convex and contains an inner point.
   For every body $C \subset \E^d$ let $C^{\circ}=\{y \in
   \DR:\scalp[c,y] \le 1$ for all $c \in C\}$ denote the
   \emph{polar} of $C$.  

   By $\LS_{j,d}$ and $\AS_{j,d}$ we denote the set of all
   $j$-dimensional linear subspaces and all $j$-dimensional affine subspaces of
   $\E^d$, respectively. For any $F \in \LS_{j,d}$ let $F^{\bot} \in
   \LS_{d-j,d}$  be the orthogonal space of $F$.  Let
   $\spann \{s_1,\dots,s_j\}$ denote 
   the \emph{linear span} $\{x \in \DR^d: x= \sum_{k=1}^j \lambda_k
   s_k, \lambda_k \in \DR\}$ of $s_1,\dots,s_j \in \DS$.
   For any set $A \in \E^d$, $\op[A]{F}$ denotes the \emph{(orthogonal)
   projection} of $A$ onto $F \in \LS_{j,d}$. If $s_1,\dots,s_j$ is an
   orthonormal basis of $F$ we also use $A_{s_1,\dots,s_j}$ instead of
   $\op[A]{F}$ and $A^{s_1,\dots,s_j}$ for the projection
   of $A$ onto $F^{\bot}$. For any $x \in \E^{d_1}$ and $y \in
   E^{d_2}$ let $x \otimes y$ denote the matrix with elements $x_i
   y_j$, $i = 1,\dots,d_1$ and $j = 1,\dots,d_2$ and note 
   that for any set of orthonormal vectors $\{s_1,\dots,s_j\}$ the
   (orthogonal) projection $P$ of $\E^d$ onto $\spann
   \{s_1,\dots,s_j\}$ can 
   be represented by the matrix $\sum_{l=1}^j s_l \otimes s_l$. 
    
   For any two sets $A,B \subseteq \E^d$ the \emph{Minkowski sum} $A+B$ is
   defined as $A+B=\{a+b \in \E^d: a \in A, b \in B\}$.
   Now for any $j \in \{1,\dots,d\}$ the \emph{inner $j$-radius}
   $r_j(C)$ of a convex set $C$ is defined by  
   $$ r_j(C)= \max \left\{\rho \ge 0 : (q+\rho \DB)\cap F \subset C, q
   \in F \in \AS_{j,d}\right\} $$
   and the \emph{outer $j$-radius $R_j(C)$} by 
   \[ R_j(C)= \min \left\{\rho \ge 0 : E + \rho \DB \supset C, E \in
   \AS_{d-j,d} \right\}.\]
   
   If a body $C_1$ arises from $C_2$ by rotation, translation and
   dilatation, we say $C_1$ is \emph{similar} to $C_2$.  Note that the
   radii of a body do not change if the body is translated or rotated;
   neither are the relationships of the radii affected by scaling the
   body.  For this reason, we will often use the word `ball' to
   signify any similar copy of $\DB$, and the same we do for
   simplices, cross-polytopes and cubes.

   Let $T^d$ denote the regular $d$-simplex
   $\conv\{e_1,\dots,e_{d+1}\}$ (embedded in $\E^{d+1}$) and
   $X^d$ the cross-polytope $\conv\{\pm e_1,\dots, \pm e_d\}$, where
   $e_k$ denotes the $k$-th unit vector (of the appropriate space). 
   By $B_{a_1,\dots,a_d}$ we denote a $d$-dimensional box of the form 
   $\left\{ x \in \DR^d: -a_i \le x_i \le a_i, i \in \{1,\dots,d\}
   \right\}$ and the cube $B_{1,\dots,1}$ we denote by  
   $B^d$.
 
\section{Regular simplices}   
   
   The following result about the inradii of simplices is taken from
   \cite{ball-92}:

   \begin{Proposition} \label{simplex_inradius}
    For the inner radii of the regular simplex of edge length
    $\sqrt{2}$ it holds
    $r_j(T^d)=\sqrt{\frac{1}{j(j+1)}}$.
   \end{Proposition}

   So we can concentrate on the outer-radii. A proof of the following
   proposition can be found in \cite{gritzmann-klee-92}:
  
   \begin{Proposition} \label{duality}
    If $C$ is a symmetric body and $1 \le j \le d$ then $r_j(C)R_j(C^0)=1$ and
    $R_j(C)r_j(C^0)=1$.
   \end{Proposition}
 
   Now we state the so far known results about outer radii of regular 
   simplices, which are taken from \cite{jung-01}, \cite{steinhagen-21}, and
   \cite{weissbach-83a} respectively.  
 
  \begin{Proposition} \label{regSimplex} \mbox{}
   \begin{enumerate}[(i)]
    \item \label{regSimplex3} $R_d(T^d) = \sqrt{\frac{d}{d+1}}$
    \item \label{regSimplex2} $R_1(T^d) = \begin{cases}
         \sqrt{\frac{1}{d+1}},& \text{if } d \text{ odd}\\
         \sqrt{\frac{d+1}{d(d+2)}},& \text{if } d \text{ even.}
        \end{cases}$  
    \item \label{regSimplex4} $R_{d-1}(T^d) =
         \sqrt{\frac{d-1}{d+1}}$ , if $d$ is odd.\footnote{Note that
           the proof of the even-case result provided in
           \cite{weissbach-83b} is wrong (see
           \cite{brandenberg-theobald-02}), but it is resettled in
           recent work \cite{brandenberg-theobald-03}}
   \end{enumerate}
  \end{Proposition}

   Proposition \ref{regSimplex} is not as complete a result as
   Proposition \ref{simplex_inradius}. At the end of this section we
   will be able to give a result on the outer radii of regular
   simplices which is much more general than the above Proposition. To
   do so, we make use of the following definition:
   
   \begin{Definition} \label{gsb-def}
    We call any set of orthonormal vectors $\{s_1,\dots,s_j\}$, $j
    \in \{1,\dots,d\}$ in $\E^{d+1}$ 
    \begin{enumerate}[(i)]
     \item a \emph{valid subset basis} (\emph{vsb} for short) if
      $\sum_{k=1}^{d+1}s_{lk}=0$ for all $l \in \{1,\dots,j\}$, and
     \item a \emph{good subset basis} (\emph{gsb} for short) if
      it is a vsb and $\sum_{l=1}^{j}s_{lk}^2=\frac{j}{d+1}$ for all $k
      \in \{1,\dots,d+1\}$.
    \end{enumerate}
   \end{Definition}

   Note that any set of orthonormal vectors $\{s_1,\dots,s_j\}$ is
   called a vsb if it spans a $j$-dimensional subspace of $\E^{d+1}_0
   = \left\{x \in \E^{d+1} : \sum_{k=1}^{d+1} x_k = 0 \right\},$
   the
   $d$-dimensional linear subspace of $\E^{d+1}$ parallel to the
   hyperplane in which we have embedded $T^d$.
   
   The projection of $T^d$ onto $\E^{d+1}_0$ can be written as
   $I^{d+1}-\frac{1}{d+1}\mathbf{1}^{d+1},$ where $I^{d+1}$ denotes
   the identity matrix in $\E^{(d+1) \times (d+1)}$ and
   $\mathbf{1}^{d+1}$ the matrix in $\E^{(d+1) \times (d+1)}$
   consisting only of 1's. Hence it holds that $\sum_{l=1}^{d} s_l
   \otimes s_l = I^{d+1} - \frac{1}{d+1}\mathbf{1}^{d+1},$
   for every
   vsb of $d$ elements. This enables us to obtain the important fact
   that each vsb is a gsb if $j=d$, which we use in Corollary
   \ref{orth-proj}.
   
   Now we start improving the results on the outer radii of regular
   simplices by giving a general lower bound, which we will prove to
   be tight in many cases further on. This theorem will also show the
   reason why we call a vsb good if it fulfills the condition (ii) in
   Definition \ref{gsb-def}.
   
  \begin{Satz} \label{weissb-gen1}
   $R_j(T^d) \ge \sqrt{\frac{j}{d+1}}$ for all $j \in
   \{1,\dots,d\}$ and equality holds iff there exists a gsb
   $\{s_1,\dots,s_j\}$ in $\E^{d+1}$. 
  \end{Satz}
  \begin{proof}
   Let $P$ denote the projection onto some subspace
   spanned by a vsb $\{s_1,\dots,s_j\}$. It follows  
   $$ \norm[Pe_k]^2 = \scalp[Pe_k,e_k] = \Scalp[\sum_{l=1}^j
   s_{lk}s_l,e_k]=\sum_{l=1}^j s_{lk}^2.$$ 
   Now assume there exists any $x \in \E^{d+1}$ such that
   $\norm[x-Pe_k]^2 < \sqrt{\frac{j}{d+1}}$ for all
   $k=1,\dots,d+1$. Summing over the $k$'s  it follows 
   \begin{align*}
    j &> \sum_{k=1}^{d+1} \norm[x-Pe_k]^2 \\
    &= \sum_{k=1}^{d+1}(\norm[x]^2-2\scalp[x,Pe_k]+\norm[Pe_k]^2)\\
    &= (d+1) \norm[x]^2 - 2 \Scalp[x,\sum_{k=1}^{d+1}\sum_{l=1}^j
       s_{lk}s_l] + \sum_{k=1}^{d+1}\sum_{l=1}^j s_{lk}^2\\
      \intertext{and since $\sum_{k=1}^{d+1} s_{lk} = 0$ and
        $\sum_{k=1}^{d+1} s_{lk}^2 = 1$}
    &= (d+1) \norm[x]^2 + j\\
    &\ge j   
   \end{align*}
   which is a contradiction. This proves the first part of the
    theorem. To prove the other part, look at the expression above; it
    is easy to see that equality in $\norm[x-Pe_k]^2 \le
    \frac{j}{d+1}$ for all $k$ can only be obtained if $x=0$ and
    $\sum_{l=1}^{j}s_{lk}^2=\frac{j}{d+1}$.
  \end{proof}

  As every vsb of $d$ vectors is already a gsb we receive the
  following corollary from Theorem \ref{weissb-gen1} from
  the basis extension property (used on $\E^{d+1}_0$):

  \begin{Folgerung} \label{orth-proj}
   For any dimension $d$ and any $j \in \{1,\dots,d-1\}$ it holds
   $R_j(T^d)=\sqrt{\frac{j}{d+1}}$ iff
   $R_{d-j}(T^d)=\sqrt{\frac{d-j}{d+1}}$ holds. Moreover the optimal
   projections take place in orthogonal subspaces.
  \end{Folgerung}
 
   Corollary \ref{orth-proj} shows that Proposition \ref{regSimplex}
   (\ref{regSimplex2}) and (\ref{regSimplex4}) correspond to each
   other in the sense that the lower bound of Theorem
   \ref{weissb-gen1} is attained in both cases for odd dimensions and
   that the bound is not attained in even dimensions.
 
  The following Proposition is a polar-version of a theorem due to
  John \cite{john-48}: 

  \begin{Proposition} \label{john}
   $\DB$ is the ellipsoid of minimal volume
   containing some body $C \subset \E^d$ iff $C \subset \DB$ and for
   some $m \ge d$ there are unit vectors $u_1,\dots,u_m$ on the
   boundary of $C$, and positive numbers $c_1,\dots,c_m$ summing to
   $d$ such that
   \begin{enumerate}[(i)]
    \item $\sum_{i=1}^m c_i u_i = 0$, and\\
    \item $\sum_{i=1}^m c_i u_i \otimes u_i = I^d$.
   \end{enumerate}
  \end{Proposition}

  It is obvious that if $C$ is a regular polytope all
  $c_i$ can be chosen as $\frac{d}{m}$ were $m$ is the number of
  vertices of $C$. But it is not obvious which other polytopes fulfill
  this property. Nevertheless, according to
   \cite{giannopoulos-milman-00} these polytopes are in an isotropic
   position, corresponding to the measure $\mu^*$ on $\DS$ that gives
   mass $\frac{d}{m}$ to all vertices $u_i$. This is the source for
   the following definition:

   \begin{Definition} \label{iso-def}
     Let $C=\conv\{u_1,\dots, u_{m}\} \subset \DB$ be a polytope,
     where all $u_i$'s are situated on $\DS$. We call $C$ \emph{quasi
       isotropic}, if all the $c_i$'s in Proposition \ref{john} can be
     taken as $\frac{d}{m}$, and \emph{isotropic}, if additionally
     $u_{i_1} \neq u_{i_2}$ for all $i_1 \neq i_2$.
   \end{Definition}  

  \begin{Lemma} \label{quasi}
   There exists a gsb $s_1,\dots,s_j$ of $\E^{d+1}$ iff there exists a
   quasi-isotropic polytope $C=\conv\{u_1,\dots, u_{d+1}\} \subset
   \E^j$, $j \le d$. Moreover if we project $T^d$ onto $\spann
   \{s_1,\dots,s_j\}$ the projection will be similar to the
   corresponding $C$.
  \end{Lemma}

  \begin{proof}
     If $C = \conv\{u_1,\dots, u_{d+1}\}$ is a quasi isotropic
     polytope then
    \begin{enumerate}[(i)]
    \item $\norm[u_k]=1$,
    \item $\sum_{k=1}^{d+1} u_k = 0$, and
    \item $\sum_{k=1}^{d+1} u_k \otimes u_k = \frac{d+1}{j} I^j$.
    \end{enumerate}
    Now let $s_l=\sqfr{j}{d+1}(u_{1,l},\dots,u_{d+1,l})^T$,
    $l=1,\dots,j$. This defines a gsb. For showing this it is
    necessary that the $s_l$ form an orthonormal set, but this is the
    case because of (iii). $\sum_{k=1}^{d+1} s_{lk}$ has to be 0, but
    this follows from (ii), and finally we need $\sum_{l=1}^j s_{lk}^2
    = \frac{j}{d+1}$ for all $k$, but this is true because of (i). The
    other direction can be shown using a similar reasoning.
    
    Now, if we project the vertices of $T^d$ onto $\spann
    \{s_1,\dots,s_j\}$ we get $P e_k = \sum_{l=1}^j s_{lk} s_l =
    \sum_{l=1}^j \sqfr{j}{d+1} u_{kl} s_l$. Hence the values
    $\sqfr{j}{d+1} u_{kl}$ are just the coordinates of the vertices of
    the projection in terms of the basis $s_1,\dots,s_j$.
   \end{proof}

  Lemma \ref{quasi} can be used in 2 ways:
  \begin{enumerate}[(i)]
   \item We know that $R_j(T^d)=\sqrt{\frac{j}{d+1}}$ whenever we find
    a quasi-isotropic $j$-dimensional polytope with $d+1$ vertices and
    vice versa (therefore, due to Proposition \ref{regSimplex}
    (\ref{regSimplex4}) there cannot be quasi-isotropic polytopes with 
    $d+2$ vertices if $d$ is odd), and
   \item we know that $R_k(C) = \sqrt{\frac{k}{j}}$ for any $k \le j$
    such that the gsb $\{s_1,\dots,s_j\}$ can be split into 2 gsb's
    $\{s_1,\dots,s_k\}$ and $\{s_{k+1},\dots,s_j\}$. 
  \end{enumerate}
  We will first concentrate our attention to (i) but come back to (ii)
  later. Because every $m$-gon is a regular body in $\E^2$, and
  because a prism
  or an anti-prism of such an $m$-gon, such that all the vertices are on
  the unit sphere, is at least isotropic, we receive (always keeping
  Corollary \ref{orth-proj} in mind) that:

  \begin{Folgerung} \label{d23}
   \begin{enumerate}[(i)]
    \item $R_2(T^d)=\sqrt{\frac{2}{d+1}}$,
      $R_{d-2}(T^d)=\sqrt{\frac{d-2}{d+1}}$ for all $d \ge 2$, and
    \item $R_3(T^d)=\sqrt{\frac{3}{d+1}}$,
      $R_{d-3}(T^d)=\sqrt{\frac{3}{d+1}}$ for all odd $d \ge3$. 
   \end{enumerate}
  \end{Folgerung}

  In the following we will not mention the ($d-j$)-cases as long as we
  make no special use of them.

  \begin{Lemma} \label{add-mult} \label{factors}
   \begin{enumerate}[(i)]
    \item \label{add} Suppose $d=d_1+d_2+1$ and
      $R_j(T^{d_1})=\sqfr{j}{d_1+1}$ and 
      $R_j(T^{d_2})=\sqfr{j}{d_2+1}$. Then it also holds
      $R_j(T^d)=\sqfr{j}{d+1}$.
    \item \label{mult} Suppose $d+1=(d_1+1)(d_2+1)$ and
      $R_{j_1}(T^{d_1})=\sqfr{j_1}{d_1+1}$ and 
      $R_{j_2}(T^{d_2})=\sqfr{j_2}{d_2+1}$ .Then
      $R_j(T^d)=\sqfr{j}{d+1}$ for all $j \in \{k_1,k_2,k_1 k_2, k_1
      (k_2+1), (k_1+1) k_2, (k_1+1)(k_2+1) \}$, where $k_i \in
      \{j_i,d_i-j_i,d_i\}$, $i=1,2$. 
   \end{enumerate}
  \end{Lemma}

  \begin{proof}
   Part (\ref{add}) is quite simple in terms of the polytopes: If $C_1$
   and $C_2$ are two quasi isotropic polytopes with $d_1+1$ and $d_2+1$
   vertices, respectively, then their convex hull has $d+1$ vertices and is
   again quasi isotropic. For part (\ref{mult}) there is a bit more to
   do. Suppose 
   $\{s_1,\dots,s_{k_1}\}$ and 
   $\{t_1,\dots,t_{k_2}\}$ are gsb's in $\E^{d_1+1}$ and
   $\E^{d_2+1}$, respectively, and consider the following three sets
   in $\E^{d+1}$: 
   $$\sqrt{\frac{1}{d_2+1}}(\underbrace{s_l^T,\dots,s_l^T}_{d_2+1})^T,
    \; l=1,\dots,k_1, $$
   and 
   $$\sqrt{\frac{1}{d_1+1}} (
   \underbrace{t_{l1},\dots,{t_{l1}}}_{d_1+1}, \dots ,
   \underbrace{t_{l(d_2+1)},\dots,{t_{l(d_2+1)}}}_{d_1+1})^T, \;
   l=1,\dots,k_2, $$ 
   and 
   $$s_{l_1} \otimes t_{l_2},  \; l_1=1,\dots,k_1, \; l_2=1,\dots,k_2
   $$
   where we take the $\otimes$-matrix as a vector, column by
   column. It is easy to 
   see that all vectors in the three sets form a vsb of size $k_1+k_2+k_1
   k_2$. Now we only need that each of the three sets forms a gsb. But
   again this is obvious for the first two sets and not hard to see for
   the last one.
  \end{proof}

  One should note that the first two groups could also be obtained by
  Part (\ref{add}) applying it $k_i$-times, and that the polytope
  corresponding to the $k_1 k_2$-gsb is a homothetic of $\conv
  \{u_{i_1} \otimes v_{i_2}, \; 
  i_1=1,\dots,d_1+1, \;  i_2=1,\dots,d_2+1\}$ if
  $\conv\{u_1,\dots,u_{d_1+1}\}$ and $\conv\{v_1,\dots,v_{d_2+1}\}$
  where the polytopes corresponding to the initial gsb's.

  The polytopes one gets out of Part (\ref{mult}) stay to be much more
  'regular' as the ones obtained Part (\ref{add}).

  \begin{Beispiel}
   Suppose $d_1=1$ and $d_2=2$ in Lemma \ref{factors}, so $d=5$. If we
   forget about the normalizing factors
   gsb's for $d_1$ and $d_2$ could be
   $$\left\{\left(\begin{array}{c} 1\\-1\end{array}\right)\right\},
   \left\{\left(\begin{array}{c} 1\\-1\\0\end{array}\right),
   \left(\begin{array}{c} 1\\1\\-2\end{array}\right) \right\},  $$  
   respectively. Using the construction in Lemma \ref{factors} we get
   the following three gsb's for $d=5$:
   \begin{enumerate}[(i)]
    \item the $d_1$-gsb by putting the $\left(1 \atop -1\right)$ vector
     $d_2+1=3$ times below each other  
     $$\left\{\left(\begin{array}{c} 1\\-1\\1\\-1\\1\\-1
          \end{array}\right)\right\},$$ 
    \item the $d_2$-gsb by taking every entry in the original
     $d_2$-gsb $d_1+1$ times below each other
     $$\left\{\left(\begin{array}{c} 1\\1\\-1\\-1\\0\\0
          \end{array}\right), \left(\begin{array}{c} 1\\1\\1\\1\\-2\\-2
          \end{array}\right)\right\},$$ and
    \item the $d_1 d_2$-gsb by multiplying each vector of the $d_1$-gsb
     coordinate wise with any vector of the $d_2$-gsb
     $$\left\{\left(\begin{array}{c} 1\\-1\\-1\\1\\0\\0
          \end{array}\right), \left(\begin{array}{c} 1\\-1\\1\\-1\\-2\\2
          \end{array}\right)\right\}.$$
   \end{enumerate}
  \end{Beispiel}
 
  Now, we derive the main theorem by making use of Corollary
  \ref{orth-proj} and Lemma  \ref{add-mult}(\ref{add}):

  \begin{Satz} \label{main}
   $R_j(T^d)=\sqfr{j}{d+1}$, if 
   \begin{enumerate}[(i)]
    \item $d$ is odd, or
    \item $j$ is even and $d \neq 2j$.
   \end{enumerate}
  \end{Satz}

  \begin{proof}
   We do an inductive proof over $j$ and $d$. We know
   already that (i) and (ii)
   are true for $j=1,2,3$. So let $j \ge 4$. Now suppose
   $d<2j$. Then $d-j<j$ and therefore the statement follows inductively 
   by applying Corollary \ref{orth-proj}, because if $d$ is odd then we do not
   depend on $j$ and if $d$ is even then $d-j$ is even if $j$ is and
   $2(d-j)=d$ would mean $2j=d$. If $d>2j$ then we can apply Lemma
   \ref{add-mult} (\ref{add}) with $R_j(T^j)$ and $R_j(T^{d-j-1})$ or in
   the case that $d-j=2j$ with $R_j(T^{j+2})$ and $R_j(T^{d-j-3})$.
   $R_j(T^j)$ and $R_j(T^{j+2})$ belong to the ($d<2j$)-case and the
   other two we can use by induction if $d-j-1$ or $d-j-3$ are good for
   one of the two cases. But if $j$ is odd we can assume that $d$ is odd and
   then this two numbers are also odd and we fulfill case (i). On the
   other hand, if $j$ is even at least one of them is not equal to $2j$
   and we have case (ii) for at least one of the two pairs. We are let
   with the case $d=2j$. But this can only be in case of even $d$ and
   hence we aren't in (i) or (ii).
  \end{proof}

  The 'only if'-direction in Theorem \ref{main} would not be true as one
  can find gsb's for the special case that $d+1=2j$, with even $j$ for
  many $d$:  

  \begin{Lemma} \label{d2j}
   In case of $d=2j$ and $j$ even  $R_j(T^d)=\sqfr{j}{d+1}$, if
   $d+1=(d_1+1)(d_2+1)$ with $d_1$ divides $j$ and if 
   $\frac{j}{d_1}$ is odd then $\frac{j}{d_1}-1  \neq \frac{d_1}{2}$. 
  \end{Lemma}

  \begin{proof}
   Because $d_1$ divides $j$ we can make use of Lemma \ref{add-mult}
   (\ref{mult}) with $j_1=d_1$ and $j_2 = \frac{j}{d_1}$ or $j_2 =
   \frac{j}{d_1}-1$ which ever is even. Now we only have to make sure
   that $2j_2 \neq d_2$. But from $2\frac{j}{d_1}=d_2$ follows
   $d+1=(d_1+1)(\frac{2j}{d_1}+1)=d+\frac{2j}{d_1}+d_1+1$ and
   therefore $j=-\frac{d_1^2}{2}$ which is a contradiction. And from
   $2\frac{j}{d_1}-1=d_2$ follows
   $d+1=(d_1+1)(\frac{2j}{d_1}-1)=d+\frac{2j}{d_1}-d_1-1$ and
   therefore $\frac{j}{d_1}=\frac{d_1}{2}+1$ the case excluded by the
   assumption. 
  \end{proof}

   \begin{table}[ht]
    \begin{tabular}{|c|*{16}{c}|}
     \hline $j,d$ & 1 & 2 & 3 & 4 & 5 & 6 & 7 & 8 & 9 & 10 & 11 & 12 &
     13 & 14 & 15 & 16 \\ \hline 1 & + & - & + & - & + & - & + & - & +
     & - & + & - & + & - & + & - \\ 2 & & + & + & + & + & + & + & + &
     + & + & + & + & + & + & + & + \\ 3 & & & + & - & + & (-) & + &
     (-) & + & (-) & + & (-) & + & (-) & + & (-) \\ 4 & & & & + & + &
     + & + & + & + & + & + & + & + & + & + & + \\ 5 & & & & & + & - &
     + & (-) & + & (-) & + & (-) & + & (-) & + & (-) \\ 6 & & & & & &
     + & + & + & + & + & + & ? & + & + & + & +\\ 7 & & & & & & & + & -
     & + & (-) & + & (-) & + & (-) & + & (-) \\ 8 & & & & & & & & + &
     + & + & + & + & + & + & + & ? \\ 9 & & & & & & & & & + & - & + &
     (-) & + & (-) & + & (-) \\ 10 & & & & & & & & & & + & + & + & + &
     + & + & + \\ \hline
    \end{tabular}
    \vspace{0.1cm}
    \caption[The existence of
      $j$-dimensional quasi isotropic polytopes with $d+1$
      vertices.]{ \label{+-?} The table shows the existence of
      $j$-dimensional quasi isotropic polytopes with $d+1$
      vertices. The first column states the $j$ value, the first row
      the value of $d$.
      A `+' indicates the existence, a `-' the non-existence. The `(-)'
      entries show that the nonexistence is not proven but very
      unlikely, the `?'-s show the open cases for even $j$. Be careful,
      in terms of the outer $j$-radii both `+' and `-', indicate that
      the radii of the regular simplices are known, and each `(-)' or
      `?' entry stands for an unsolved case.}
   \end{table}

  Lemma \ref{d2j} includes the case that 3 divides $d+1$ (because $j$
  is even) and the case
  that $d+1=(d_1+1)^2$ with 4 does not divide $d_1$ (because
  $j=\frac{d_1}{2}+1$). 

  On the other hand Lemma \ref{add-mult} (\ref{mult}) cannot help in the case
  $d=2j$, $j$ even if $d+1$ is prime neither does it help always if
  $d+1$ is not prime since if $d+1=5*17$ we have $j=42$. Here we
  would need $j_1 \in \{2,4\}$, but if $j_1=2$ it follows $j_2 > 16$
  and $j_1=4$ is not possible because neither 4 nor 5 divides 42.

\section{Boxes and cross-polytopes}

   \begin{Proposition} \label{box}
     Let $0< a_1 \le \dots \le a_d$. Then
    \begin{enumerate}[(i)] 
    \item
      $$r_j(B_{a_1,\dots,a_d}) = \sqrt{\frac{a_1^2 + \dots +
          a_{d-k}^2}{j-k}},$$
      where $k$ is the smallest of the
      integers $0,\dots,j-1$ that satisfies $$a_{d-k} \le
      \sqrt{\frac{a_1^2 + \dots + a_{d-k-1}^2}{j-k-1}},$$
      and
    \item $$R_j(X_{a_1,\dots,a_d}) = \sqfr{(j-k)\prod_{i=k}^d
        a_i^2}{\sum_{i=k}^d \prod_{l \neq i} a_l^2},$$
      where $k$ is
      the smallest of the integers $0,\dots,j-1$ that satisfies
      $$a_k \ge \sqfr{(j-k-1)\prod_{i=k+1}^d a_i^2}{\sum_{i=k+1}^d
        \prod_{l \neq i} a_l^2}.$$
    \end{enumerate}
   \end{Proposition}
   
   The corresponding result about the outer radii of boxes seems to be
   very intuitive. It says that one should just project the box
   through one of its smallest faces. One gets the inner radii of
   cross-polytopes from polarization. Before we state the final
   theorem, we give a technical lemma which will be useful in the
   proof of the theorem.

   \begin{Lemma} \label{pm-config}
     Let $s_l \in \E^d$, $l=1,\dots,j$, $j \le d$, $d \ge 2$ be a set
     of orthonormal vectors and $a_1,\dots,a_d \in \DR_+$. Then there
     exists a choice of plus and minus signs in
     $\sum_{l=1}^j(\sum_{k=1}^d \pm a_k s_{lk})^2$ such that this is
     at least $\sum_{k=1}^d a_k^2 \sum_{l=1}^j s_{lk}^2$.
   \end{Lemma}
   
   \begin{proof}
     Essentially we have to show that there is an
     $\alpha=(\alpha_1,\dots,\alpha_d)$ with $\alpha_k \in
     \{-a_k,a_k\}$ such that $$\Gamma_\alpha := \sum_{1 \le k_1 < k_2
       \le d} \alpha_{k_1} \alpha_{k_2} \zeta_{k_1,k_2} \ge 0, \text{
       where } \zeta_{k_1,k_2}:=\sum_{l=1}^j s_{lk_1}s_{lk_2}.$$
     This
     is done by an inductive proof.
     
     First consider the case $d=2$. Then $\Gamma_\alpha= \alpha_1
     \alpha_2 \zeta_{1,2}$. So we choose $\alpha_i=a_i$, $i=1,2$ if
     $\zeta_{1,2} \ge 0$ and if $\zeta_{1,2} < 0$ we choose $\alpha_1
     = a_1$ and $\alpha_2 = - a_2$.
    
     Now let $d=3$. Hence $\Gamma_\alpha= \alpha_1 \alpha_2
     \zeta_{1,2} + \alpha_1 \alpha_3 \zeta_{1,3} + \alpha_2 \alpha_3
     \zeta_{2,3}$.  Now suppose $\Gamma_{(a_1,a_2,-a_3)} < 0$ and
     $\Gamma_{(-a_1,a_2,-a_3)} < 0$. It follows that $$0 >
     \Gamma_{(a_1,a_2,-a_3)} + \Gamma_{(-a_1,a_2,-a_3)} = -2a_2 a_3
     \zeta_{2,3}$$
     and therefore that $\zeta_{2,3} > 0$. Analogously
     one can show that $\zeta_{1,2}$ and $\zeta_{1,3}$ are positive;
     but then we can choose $\alpha=(a_1,a_2,a_3)$.
     
     By knowing that the statement is correct for $d=2,3$ we can take
     an inductive step of 2, that means we assume the statement is
     proven up to some $d$ and now conclude that it is also true for
     $d+2$.
    
     Now suppose the statement would be wrong for $d+2$, meaning
     $\Gamma_\alpha < 0$ for all possible choices of $\alpha \in
     \E^{d+2}$. Hence
    \begin{align*}
      0 & > \Gamma_{\alpha_1,\dots,\alpha_d,a_{d+1},a_{d+2}} +
      \Gamma_{\alpha_1,\dots,\alpha_d,-a_{d+1},a_{d+2}} +
      \Gamma_{\alpha_1,\dots,\alpha_d,a_{d+1},-a_{d+2}} +
      \Gamma_{\alpha_1,\dots,\alpha_d,-a_{d+1},-a_{d+2}} \\
      & = 4 \sum_{1 \le k_1 < k_2 \le d} \alpha_{k_1} \alpha_{k_2}
      \zeta_{k_1,k_2}.
    \end{align*}
    However, this is not possible as by the induction hypothesis
    $$\sum_{1 \le k_1 < k_2 \le d} \alpha_{k_1} \alpha_{k_2}
    \zeta_{k_1,k_2} \ge 0$$
    for at least one possible choice of
    $\alpha$.
   \end{proof}

   \begin{Satz} \label{box2}
     Let $0< a_1 \le \dots \le a_d$. Then
    \begin{enumerate}[(i)]
    \item $R_j(B_{a_1,\dots,a_d})=\sqrt{a_1^2+\dots+a_j^2}$, and
    \item $r_j(X_{a_1,\dots,a_d})=\frac{\prod_{i=d-j+1}^d
        a_i}{\sqrt{\sum_{i=d-j+1}^d \prod_{l \neq i} a_l^2}}$.
    \end{enumerate}
   \end{Satz}
   
   \begin{proof}
     It suffices to show Part (i), Part (ii) follows then from
     Proposition \ref{duality}, and as the result is obvious if $d=1$
     we can assume that $d \ge 2$.  Any vertex $v$ of
     $B_{a_1,\dots,a_d}$ can be written in the form $v=\sum_{k=1}^d
     \pm a_k e_k$ and all possible choices of the plus and minuses in
     that formular leads to a vertex of $B_{a_1,\dots,a_d}$. Hence,
     for every projection $P=\sum_{l=1}^j s_l \otimes s_l$ with
     pairwise orthogonal unit-vectors $s_l \in \E^d$, it holds that
     $\norm[Pv]^2 = \sum_{l=1}^j \scalp[v,s_l]^2 =
     \sum_{l=1}^j(\sum_{k=1}^d \pm a_k s_{lk})^2 $ and because of
     Lemma \ref{pm-config} there exists a vertex of
     $B_{a_1,\dots,a_d}$ such that this is at least $\sum_{k=1}^d
     a_k^2 \sum_{l=1}^j s_{lk}^2$. Now extend the set
     $\{s_1,\dots,s_j\}$ to an orthonormal basis of $\E^d$. As
     $\sum_{l=1}^d s_l \otimes s_l=I$ it follows that $\sum_{k=1}^d
     s_{lk}^2 = \sum_{l=1}^d s_{lk}^2 = 1$, for all $k=1,\dots,d$, and
     therefore that $t_k := \sum_{l=1}^j s_{lk}^2 \in [0,1]$. Now,
     because $\sum_{k=1}^d t_k = \sum_{l=1}^j \sum_{k=1}^d s_{lk}^2$
     has to equal $j$ the minimum value of $\sum_{k=1}^d t_k a_k^2$
     will be achieved for $t_1 = \dots = t_j = 1$ and
     $t_{j+1}=\dots=t_d=0$. Hence $R_j(B_{a_1,\dots,a_d}) \ge
     \sqrt{a_1^2+\dots+a_j^2}$. But as the projection of
     $B_{a_1,\dots,a_d}$ through its $j$-face $B_{a_1,\dots,a_j}$
     achieves this value we get the desired result.
   \end{proof}   
   
   Of course, one can easily get the radii of cubes and regular
   cross-polytopes from the results about the radii of general boxes
   and cross-polytopes:
  
   \begin{Folgerung} \label{outer-cross-polytope}
     The following hold:
    \begin{enumerate}[(i)]
    \item $r_j(B^d)= \sqrt{\frac{d}{j}}$, and
    \item $R_j(X^d)=\sqrt{\frac{j}{d}}$.
    \item $R_j(B^d)=\sqrt{j}$, and
    \item $r_j(X^d)=\sqfr{1}{j}$.
    \end{enumerate}
   \end{Folgerung}

   \begin{proof}
     Part (i) and (ii) follow from Proposition \ref{box} by choosing
     $k=0$ there, Part (iii) and (iv) from Theorem \ref{box2}.
   \end{proof}
   
   However, one should recognize that we can prove Corollary
   \ref{outer-cross-polytope} almost without using Proposition
   \ref{box} and \ref{box2}.  Except in the case where $2j=d-1$ and
   Lemma \ref{d2j} does not hold we can use this lemma and Theorem
   \ref{main} to show (ii) and Theorem \ref{main} suffices to show
   (iii). Parts (i) and (iv) would follow again by duality.
   
   How do we do this in detail? Part (iii) follows from the fact that
   the cube and all its faces (which are again cubes) are quasi
   isotropic, that by projecting through a face of a cube all vertices
   stay on the circumsphere, and that the distance from the center of
   the cube to the center of any of its $j$-faces is $\sqrt{d-j}$.
   
   To prove Part (ii) we remember the second statement after Lemma
   \ref{quasi}. First, we project $T^{2d-1}$ onto $\sqfr{1}{2} X^d$ by
   using the gsb $$\sqfr{1}{2}\left( {s_1 \atop -s_1}
   \right),\dots,\sqfr{1}{2}\left({s_{d-1} \atop -s_{d-1}} \right),
   \left({\DEins_{d-1} \atop - \DEins_{d-1}}\right) \ ,$$
   where
   $s_1,\dots,s_{d-1}$ is any gsb for $T^{d-1}$. Now because for every
   even $j$, which is not excluded by both the theorem and the lemma,
   there exists a subset of size $j$ of $s_1,\dots,s_{d-1}$, wlog
   $s_1, \dots, s_j$. But hence the sets $$\sqfr{1}{2}\left( {s_1
       \atop -s_1} \right),\dots,\sqfr{1}{2}\left({s_j \atop -s_j}
   \right)$$
   and
   $$\sqfr{1}{2}\left( {s_1 \atop -s_1} \right), \dots
   ,\sqfr{1}{2}\left({s_j \atop -s_j} \right), \left({\DEins_{d-1}
       \atop -\DEins_{d-1}}\right)$$
   are gsb's in $\E^{2d}$ and
   therefore there exists a projection of $X^d$ onto any $j'$ subspace
   such that it attains the lower bound $\sqfr{j'}{d}$, except the
   case where $2j'$ or $2j'-1$ does not pass the conditions of Theorem
   \ref{main} or Lemma \ref{d2j}.

\bibliography{references-after-sub}
\bibliographystyle{plain}

\end{document}